\documentclass[11pt]{amsart}
\usepackage{amssymb}
\usepackage{amsmath}
\usepackage[active]{srcltx}
\usepackage{t1enc}
\usepackage[latin2]{inputenc}
\usepackage{verbatim}
\usepackage{amsmath,amsfonts,amssymb,amsthm}
\usepackage[mathcal]{eucal}
\usepackage{enumerate}
\usepackage[centertags]{amsmath}
\usepackage{graphics}

\newtheorem{theorem}{Theorem}

\begin{document}
\author{George Tephnadze}
\title[Vilenkin-Fourier coefficients]{A note on the Vilenkin-Fourier
coefficients }
\address{G. Tephnadze, Department of Mathematics, Faculty of Exact and
Natural Sciences, Tbilisi State University, Chavchavadze str. 1, Tbilisi
0128, Georgia}
\email{giorgitephnadze@gmail.com}
\date{}
\maketitle

\begin{abstract}
The main aim of this paper is to find the estimation for Vilenkin-Fourier
coefficients.
\end{abstract}

\date{}



Let $P_{+}$ denote the set of the positive integers, $P:=P_{+}\cup \{0\}.$

Let $m:=(m_{0,}m_{1....})$ denote a sequence of the positive integers not
less than 2.

Denote by
\begin{equation*}
Z_{m_{k}}:=\{0,1,...m_{k}-1\}
\end{equation*}
the additive group of integers modulo $m_{k}.$

Define the group $G_{m}$ as the complete direct product of the group $%
Z_{m_{j}}$ with the product of the discrete topologies of $Z_{m_{j}}$ $^{,}$%
s.

The direct product $\mu $ of the measures
\begin{equation*}
\mu _{k}\left( \{j\}\right) :=1/m_{k},\text{ \qquad }(j\in Z_{m_{k}})
\end{equation*}%
is the Haar measure on $G_{m_{\text{ }}}$with $\mu \left( G_{m}\right) =1.$

If $\sup_{n}m_{n}<\infty $, then we call $G_{m}$ a bounded Vilenkin group.
If the generating sequence $m$ is not bounded then $G_{m}$ is said to be an
unbounded Vilenkin group. \textbf{In this paper we discuss bounded Vilenkin
groups only.}

The elements of $G_{m}$ represented by sequences
\begin{equation*}
x:=(x_{0},x_{1,...,}x_{j,...}),\qquad \left( \text{ }x_{k}\in
Z_{m_{k}}\right) .
\end{equation*}

It is easy to give a base for the neighborhood of $G_{m}$
\begin{equation*}
I_{0}\left( x\right) :=G_{m},
\end{equation*}%
\begin{equation*}
I_{n}(x):=\{y\in G_{m}\mid y_{0}=x_{0},...y_{n-1}=x_{n-1}\},\text{ \ }(x\in
G_{m},\text{ }n\in P).
\end{equation*}%
Denote $I_{n}:=I_{n}\left( 0\right) $ for $n\in P$ and $\overset{-}{I_{n}}%
:=G_{m}$ $\backslash $ $I_{n}$ .

It is evident

\begin{equation}
\overset{-}{I_{N}}=\overset{N-1}{\underset{s=0}{\bigcup }}I_{s}\backslash
I_{s+1}.  \label{2}
\end{equation}

If we define the so-called generalized number system based on $m$ in the
following way :
\begin{equation*}
M_{0}:=1,\text{ \qquad }M_{k+1}:=m_{k}M_{k\text{ }},\ \qquad (k\in P),
\end{equation*}%
then every $n\in P$ can be uniquely expressed as $n=\overset{\infty }{%
\underset{k=0}{\sum }}n_{j}M_{j},$ where $n_{j}\in Z_{m_{j}}$ $~(j\in P)$
and only a finite number of $n_{j}`$s differ from zero. Let $\left\vert
n\right\vert :=\max $ $\{j\in P;$ $n_{j}\neq 0\}.$

Denote by $L_{1}\left( G_{m}\right) $ the usual (one dimensional) Lebesque
space.

Next, we introduce on $G_{m}$ an ortonormal system which is called the
Vilenkin system.

At first define the complex valued function $r_{k}\left( x\right)
:G_{m}\rightarrow C,$ the generalized Rademacher functions as
\begin{equation*}
r_{k}\left( x\right) :=\exp \left( 2\pi ix_{k}/m_{k}\right) ,\text{ \qquad }%
\left( i^{2}=-1,\text{ }x\in G_{m},\text{ }k\in P\right) .
\end{equation*}

Now define the Vilenkin system $\psi :=(\psi _{n}:n\in P)$ on $G_{m}$ as:
\begin{equation*}
\psi _{n}(x):=\overset{\infty }{\underset{k=0}{\Pi }}r_{k}^{n_{k}}\left(
x\right) ,\text{ \qquad }\left( n\in P\right) .
\end{equation*}

Specifically, we call this system the Walsh-Paley one if m=2.

The Vilenkin system is ortonormal and complete in $L_{2}\left( G_{m}\right)
\,$\cite{AVD,Vi}.

Now we introduce analogues of the usual definitions in Fourier-analysis. If $%
f\in L_{1}\left( G_{m}\right) $ we can establish the the Fourier
coefficients, the partial sums, the Dirichlet kernels with respect to the
Vilenkin system in the usual manner:
\begin{eqnarray*}
\widehat{f}\left( k\right) &:&=\int_{G_{m}}f\overline{\psi }_{k}d\mu ,\text{%
\thinspace }\left( \text{ }k\in P\text{ }\right) , \\
S_{n}f &:&=\sum_{k=0}^{n-1}\widehat{f}\left( k\right) \psi _{k},\left( \text{
}n\in P_{+},\text{ }S_{0}f:=0\right) , \\
D_{n} &:&=\sum_{k=0}^{n-1}\psi _{n},\left( \text{ }n\in P_{+}\text{ }\right)
.
\end{eqnarray*}

Recall that (see \cite{AVD})
\begin{equation}
\quad \hspace*{0in}D_{M_{n}}\left( x\right) =\left\{
\begin{array}{l}
\text{ }M_{n}\text{\thinspace ,\thinspace \thinspace \thinspace if\thinspace
\thinspace }x\in I_{n}, \\
\text{ }0\text{\thinspace ,\thinspace \thinspace \thinspace \thinspace if
\thinspace \thinspace }x\notin I_{n}.%
\end{array}%
\right.   \label{3}
\end{equation}%
\vspace{0pt}and
\begin{equation}
D_{n}\left( x\right) =\psi _{n}\left( x\right) \left( \sum_{j=0}^{\infty
}D_{M_{j}}\left( x\right) \sum_{u=m_{j}-n_{j}}^{m_{j}-1}r_{j}^{u}\left(
x\right) \right) .  \label{3aa}
\end{equation}%
The norm (or quasinorm) of the space $L_{p}(G_{m})$ is defined by \qquad

\begin{equation*}
\left\Vert f\right\Vert _{p}:=\left( \int_{G_{m}}\left\vert f\right\vert
^{p}d\mu \right) ^{1/p}\qquad \left( 0<p<\infty \right) .
\end{equation*}

The $\sigma -$algebra generated by the intervals $\left\{ I_{n}\left(
x\right) :x\in G_{m}\right\} $ will be denoted by $\digamma _{n}$ $\left(
n\in P\right) .$ The conditional expectation operators relative to $\digamma
_{n}\left( n\in P\right) $ are denoted by $E_{n}.$ Then

\begin{equation*}
E_{n}f\left( x\right) =S_{M_{n}}f\left( x\right) =\sum_{k=0}^{M_{n}-1}%
\widehat{f}\left( k\right) w_{k}=\left| I_{n}\left( x\right) \right|
^{-1}\int_{I_{n}\left( x\right) }f(x)d\mu (x),
\end{equation*}
where $\left| I_{n}\left( x\right) \right| =M_{n}^{-1}$ denotes the length
of $I_{n}\left( x\right) .$

A sequence $F=\left( f^{\left( n\right) },\text{ }n\in P\right) $ of
functions $f^{\left( n\right) }\in L_{1}\left( G\right) $ is said to be a
dyadic martingale if (for details see e.g. \cite{We1})

$\left( i\right) $ $f^{\left( n\right) }$ is $\digamma _{n}$ measurable for
all $n\in P,$

$\left( ii\right) $ $E_{n}f^{\left( m\right) }=f^{\left( n\right) }$ for all
$n\leq m.$

The maximal function of a martingale $f$ is denoted by \qquad
\begin{equation*}
f^{*}=\sup_{n\in P}\left| f^{\left( n\right) }\right| .
\end{equation*}

In case $f\in L_{1},$ the maximal functions are also be given by
\begin{equation*}
f^{*}\left( x\right) =\sup_{n\in P}\frac{1}{\left| I_{n}\left( x\right)
\right| }\left| \int_{I_{n}\left( x\right) }f\left( u\right) \mu \left(
u\right) \right| .
\end{equation*}

For $0<p<\infty $ the Hardy martingale spaces $H_{p}$ $\left( G_{m}\right) $
consist of all martingales for which
\begin{equation*}
\left\| f\right\| _{H_{p}}:=\left\| f^{*}\right\| _{L_{p}}<\infty .
\end{equation*}

If $f\in L_{1},$ then it is easy to show that the sequence $\left(
S_{M_{n}}f:n\in P\right) $ is a martingale. If $f=\left( f^{\left( n\right)
},n\in P\right) $ is martingale then the Vilenkin-Fourier coefficients must
be defined in a slightly different manner: $\qquad \qquad $
\begin{equation}
\widehat{f}\left( i\right) :=\lim_{k\rightarrow \infty
}\int_{G_{m}}f^{\left( k\right) }\left( x\right) \overline{\Psi }_{i}\left(
x\right) d\mu \left( x\right) .  \label{3a}
\end{equation}
\qquad \qquad \qquad \qquad

The Vilenkin-Fourier coefficients of $f\in L_{1}\left( G_{m}\right) $ are
the same as those of the martingale $\left( S_{M_{n}}f:n\in P\right) $
obtained from $f$ .

A bounded measurable function $a$ is p-atom, if there exist a interval $I$,
such that \qquad
\begin{equation}
\left\{
\begin{array}{l}
a)\qquad \int_{I}ad\mu =0, \\
b)\ \qquad \left\| a\right\| _{\infty }\leq \mu \left( I\right) ^{-1/p}, \\
c)\qquad \text{supp}\left( a\right) \subset I.\qquad%
\end{array}
\right.  \label{8a}
\end{equation}

The Hardy martingale spaces $H_{p}$ $\left( G\right) $ for $0<p\leq 1$ have
an atomic atomic characterization (see \cite{We1}):

\textbf{Theorem W 1. }A martingale $f=\left( f^{\left( n\right) },n\in
P\right) $ is in $H_{p}\left( 0<p\leq 1\right) $ if and only if there exist
a sequence $\left( a_{k},k\in P\right) $ of p-atoms and a sequence $\left(
\mu _{k},k\in P\right) $ of a real numbers such that for every $n\in P$

\begin{equation}
\qquad \sum_{k=0}^{\infty }\mu _{k}S_{M_{n}}a_{k}=f^{\left( n\right) },
\label{11aa}
\end{equation}

\begin{equation*}
\qquad \sum_{k=0}^{\infty }\left| \mu _{k}\right| ^{p}<\infty .
\end{equation*}
Moreover, $\left\| f\right\| _{H_{p}}\backsim \inf \left( \sum_{k=0}^{\infty
}\left| \mu _{k}\right| ^{p}\right) ^{1/p}$, where the infimum is taken over
all decomposition of $f$ of the form (\ref{11aa}).

It follows that if any operator is uniformly bounded on the p-atom, then it
is a bounded operator from the Hardy space $H_{p}$ to the space $L_{p}.$
Moreover, the following theorem is true (see \cite{We3}):

\textbf{Theorem W 2. }Suppose that an operator $T$ is sublinear and for some
$0<p\leq 1$

\begin{equation*}
\int\limits_{\overset{-}{I}}\left\vert Ta\right\vert ^{p}d\mu \leq
c_{p}<\infty ,
\end{equation*}
for every $p$-atom $a$ ,where $I$ denote the support of the atom. If $T$ is
bounded from $L_{\infty \text{ }}$ to $L_{\infty \text{ }},$then
\begin{equation*}
\left\Vert Tf\right\Vert _{L_{p}\left( G_{m}\right) }\leq c_{p}\left\Vert
f\right\Vert _{H_{p}\left( G_{m}\right) }.
\end{equation*}

The classical inequality of Hardy type is well known in the trigonometric as
well as in the Vilenkin-Fourier analysis. Namely,
\begin{equation*}
\overset{\infty }{\underset{k=1}{\sum }}\frac{\left\vert \widehat{f}\left(
k\right) \right\vert }{k}\leq c\left\Vert f\right\Vert _{H_{1}},
\end{equation*}
where the function $f$ belongs to the Hardy space $H_{1}$ and $c$ is an
absolute constant. This was proved in the trigonometric case by Hardy and
Littlewood \cite{hl} (see also Coifman and Weiss \cite{cw}) and for Walsh
system in \cite{sws}.

Weisz \cite{We1,We4} generalized this result for Vilenkin system and proved:
\begin{equation}
\overset{\infty }{\underset{k=1}{\sum }}\frac{\left| f\left( k\right)
\right| ^{p}}{k^{2-p}}\leq c\left\| f\right\| _{H_{p}}^{p},  \label{1bb}
\end{equation}
for all $f\in H_{p}$ $\left( 0<p\leq 2\right) .$

It is also well-known (see \cite{AVD}) that
\begin{equation*}
\widehat{f}\left( n\right) \rightarrow 0,\text{ when }n\rightarrow \infty ,
\end{equation*}
for all $\ f\in L_{1}\left( G_{m}\right) ,$ where $\widehat{f}\left(
n\right) $ denotes n-th Fourier coefficients of the function $f.$

The main aim of this paper is to prove that the following theorem is true:

\begin{theorem}
a) Let $0<p<1$ and $f\in H_{p}\left( G_{m}\right) .$ Then there exists an
absolute constant $c_{p}$, defend only $p,$ such that
\end{theorem}

\begin{equation*}
\left\vert \widehat{f}\left( n\right) \right\vert \leq
c_{p}n^{1/p-1}\left\Vert f\right\Vert _{H_{p}}.
\end{equation*}

b) \textbf{\ }Let $0<p<1$ and $\Phi \left( n\right) $ is any nondecreasing,
nonnegative function, satisfying condition
\begin{equation*}
\overline{\underset{n\rightarrow \infty }{\lim }}\frac{n^{1/p-1}}{\Phi
\left( n\right) }=\infty ,
\end{equation*}
then there exists a martingale $f_{0}\in H_{p}\left( G_{m}\right) ,$ such
that
\begin{equation*}
\underset{n\rightarrow \infty }{\overline{\lim }}\frac{\left| \widehat{f}%
_{0}\left( n\right) \right| }{\Phi \left( n\right) }=\infty .
\end{equation*}

\textbf{Proof of Theorem 1. }At first we prove that the maximal operator
\begin{equation*}
\overset{\sim }{S_{p}^{\ast }}\left( f\right) :=\sup_{n\in P}\frac{%
\left\vert S_{n}f\right\vert }{\left( n+1\right) ^{1/p-1}}
\end{equation*}
is bounded from the Hardy space $H_{p\text{ }}$ to the space $L_{p}$ for $%
0<p<1.$ Since the maximal operator $\overset{\sim }{S_{p}^{\ast }}$ is
bounded from $L_{\infty }$ to $L_{\infty }$ by Theorem W 2 we obtain that
the proof of theorem will be complete, if we show that

\begin{equation*}
\int\limits_{\overline{I}_{N}}\left( \underset{n\in P}{\sup }\frac{\left|
S_{n}a\right| }{\left( n+1\right) ^{1/p-1}}\right) ^{p}d\mu \leq c<\infty
\text{ },\text{\qquad when }0<p<1,
\end{equation*}
for every $p$-atom $a$ $\left( 0<p\leq 1\right) ,$ where $I$ denote the
support of the atom$.$ Let $a$ be an arbitrary $p$-atom with support$\ I$
and $\mu \left( I\right) =M_{N}^{-1}.$ We may assume that $I=I_{N}$ $.$ It
is easy to see that $S_{n}\left( a\right) =0$ when $n\leq M_{N}$ . Therefore
we can suppose that $n>M_{N}.$

Let $x\in I_{l}\backslash I_{l+1}.$ Combining (\ref{3}) and (\ref{3aa}) we
have

\begin{equation}
\left| D_{n}\left( x\right) \right| \leq \underset{j=0}{\overset{l}{\sum }}%
n_{j}D_{M_{j}}\left( x\right) =\underset{j=0}{\overset{l}{\sum }}%
n_{j}M_{j}\leq cM_{l}.  \label{11}
\end{equation}

Since $t\in I_{N}$ and $x\in I_{s}\backslash I_{s+1},$ $s=0,...N-1,$ we
obtain that $x-t\in I_{s}\backslash I_{s+1}$. Using (\ref{11}) we get

\begin{equation*}
\left| D_{n}\left( x-t\right) \right| \leq cM_{s},
\end{equation*}
and

\begin{equation}
\int_{I_{N}}\left| D_{n}\left( x-t\right) \right| d\mu \left( t\right) \leq
c\,\,\,\frac{M_{s}}{M_{N}}.  \label{11a}
\end{equation}

Hence

\begin{eqnarray}
&&\frac{\left| S_{n}\left( a\right) \right| }{\left( n+1\right) ^{1/p-1}}
\label{12} \\
&\leq &\int_{I_{N}}\left| a\left( t\right) \right| \left| \frac{D_{n}\left(
x-t\right) }{\left( n+1\right) ^{1/p-1}}\right| d\mu \left( t\right)  \notag
\\
&\leq &\frac{\left\| a\right\| _{\infty }}{\left( n+1\right) ^{1/p-1}}%
\int_{I_{N}}\left| D_{n}\left( x-t\right) \right| d\mu \left( t\right)
\notag \\
&\leq &\frac{M_{N}^{1/p}}{\left( n+1\right) ^{1/p-1}}\int_{I_{N}}\left|
D_{n}\left( x-t\right) \right| d\mu \left( t\right)  \notag \\
&\leq &\frac{cM_{N}^{1/p}}{M_{N}^{1/p-1}}\frac{M_{s}}{M_{N}}=cM_{s}.  \notag
\end{eqnarray}

Combining (\ref{2}) and (\ref{12}) we have
\begin{eqnarray*}
&&\int_{\overline{I_{N}}}\left| \widetilde{S}_{p}^{*}a\left( x\right)
\right| ^{p}d\mu \left( x\right) \\
&=&\overset{N-1}{\underset{s=0}{\sum }}\int_{I_{s}\backslash I_{s+1}}\left|
\widetilde{S}_{p}^{*}a\left( x\right) \right| ^{p}d\mu \left( x\right) \\
&\leq &c\overset{N-1}{\underset{s=0}{\sum }}\frac{M_{s}^{p}}{M_{s}}<c<\infty
.
\end{eqnarray*}

Now we are ready to prove the main result.\textbf{\ }Let $0<p<1.$Then

\begin{equation*}
\left\vert \widehat{f}\left( n\right) \right\vert =\left\vert
S_{n+1}f-S_{n}f\right\vert \leq 2\underset{n\in P}{\sup }\left\vert
S_{n}f\right\vert
\end{equation*}
and
\begin{equation}
\frac{\left\vert \widehat{f}\left( n\right) \right\vert }{\left( n+1\right)
^{1/p-1}}\leq 2\underset{n\in P}{\sup }\frac{\left\vert S_{n}f\right\vert }{%
\left( n+1\right) ^{1/p-1}},  \label{15}
\end{equation}
Consequently,
\begin{equation*}
\frac{\left\vert \widehat{f}\left( n\right) \right\vert }{\left( n+1\right)
^{1/p-1}}\leq 2\left\Vert \underset{n\in P}{\sup }\frac{\left\vert
S_{n}f\right\vert }{\left( n+1\right) ^{1/p-1}}\right\Vert _{p}\leq
c_{p}\left\Vert f\right\Vert _{H_{p}}.
\end{equation*}
It follows that

\begin{equation*}
\left\vert \widehat{f}\left( n\right) \right\vert \leq
c_{p}n^{1/p-1}\left\Vert f\right\Vert _{H_{p}}.
\end{equation*}

\textbf{b) }In the proof of second part of theorem we follow the method of
\textit{\ }Blahota, Gát and Goginava (see \cite{gog1, gog2} ).\textbf{\ }

Let $0<p<1$ and $\Phi \left( n\right) $ is any nondecreasing, nonnegative
function, satisfying condition
\begin{equation}
\overline{\underset{n\rightarrow \infty }{\lim }}\frac{n^{1/p-1}}{\Phi
\left( n\right) }=\infty ,  \label{20}
\end{equation}
then for every $0<p<1,$ there exists an increasing sequence $\left\{ \alpha
_{k}:k\in P\right\} $ of the positive integers such that:

\qquad
\begin{equation}
\sum_{k=0}^{\infty }\left( \frac{\Phi \left( M_{\alpha _{k}}\right) }{%
M_{\alpha _{k}}^{1/p-1}}\right) ^{p/2}<\infty  \label{2aa}
\end{equation}

Let \qquad
\begin{equation*}
f^{\left( A\right) }\left( x\right) =\sum_{\left\{ k;\text{ }\alpha
_{k}<A\right\} }\lambda _{k}a_{k},
\end{equation*}
where
\begin{equation*}
\lambda _{k}=\left( \frac{\Phi \left( M_{\alpha _{k}}\right) }{M_{\alpha
_{k}}^{1/p-1}}\right) ^{1/2}
\end{equation*}
and

\begin{equation*}
a_{k}\left( x\right) :=\frac{M_{\alpha _{k}}^{1/p-1}}{M}\left( D_{M_{\alpha
_{k}+1}}\left( x\right) -D_{M_{_{\alpha _{k}}}}\left( x\right) \right) .
\end{equation*}

It is easy to show that the martingale $\,f=\left( f^{\left( 1\right)
},f^{\left( 2\right) }...f^{\left( A\right) }...\right) \in H_{p}.$

Indeed, since

\begin{equation}
S_{M_{A}}a_{k}\left( x\right) =\left\{
\begin{array}{l}
a_{k}\left( x\right) ,\text{ \qquad }\alpha _{k}<A, \\
0,\text{ \qquad }\alpha _{k}\geq A,%
\end{array}
\right.  \label{4}
\end{equation}

\begin{equation*}
\text{supp}(a_{k})=I_{\alpha _{k}},
\end{equation*}

\begin{equation*}
\int_{I_{\alpha _{k}}}a_{k}d\mu =0
\end{equation*}

and

\begin{equation*}
\left\Vert a_{k}\right\Vert _{\infty }\leq \frac{M_{\alpha _{k}}^{1/p-1}}{M}%
M_{\alpha _{k}+1}\leq (M_{\alpha _{k}})^{1/p}=(\text{supp }a_{k})^{-1/p}
\end{equation*}%
if we apply Theorem W 1 and (\ref{2aa}) we conclude that $f\in H_{p}.$

It is easy to show that

\begin{equation}
\widehat{f}(j)=\left\{
\begin{array}{l}
\frac{1}{M}M_{\alpha _{k}}^{\left( 1/p-1\right) /2}\Phi ^{1/2}\left(
M_{\alpha _{k}}\right) ,\,\,\text{ if \thinspace \thinspace }j\in \left\{
M_{\alpha _{k}},...,\text{ ~}M_{\alpha _{k}+1}-1\right\} ,\text{ }k=0,1,2...
\\
0,\text{ \thinspace \thinspace \thinspace if \thinspace \thinspace
\thinspace }j\notin \bigcup\limits_{k=1}^{\infty }\left\{ M_{\alpha
_{k}},...,\text{ ~}M_{\alpha _{k}+1}-1\right\} \text{.}%
\end{array}
\right.  \label{5}
\end{equation}

It follows that
\begin{eqnarray*}
\underset{n\rightarrow \infty }{\overline{\lim }}\frac{\widehat{f}(n)}{\Phi
\left( n\right) } &\geq &\lim_{k\rightarrow \infty }\frac{\widehat{f}%
(M_{\alpha _{k}})}{\Phi \left( M_{\alpha _{k}}\right) } \\
&\geq &\underset{k\rightarrow \infty }{\lim }\frac{1}{M}\frac{M_{\alpha
_{k}}^{\left( 1/p-1\right) /2}\Phi ^{1/2}\left( M_{\alpha _{k}}\right) }{%
\Phi \left( M_{\alpha _{k}}\right) } \\
&\geq &\underset{k\rightarrow \infty }{\lim }\left( \frac{M_{\alpha
_{k}}^{\left( 1/p-1\right) }}{\Phi \left( M_{\alpha _{k}}\right) }\right)
^{1/2}=\infty .
\end{eqnarray*}

Theorem 1 is proved.

\end{document}